\documentclass{elsart3-1}

\usepackage{amssymb}
\usepackage{amsmath}
\usepackage[english,francais]{babel}

\newtheorem{theorem}{Theorem}[section]
\newtheorem{lemma}[theorem]{Lemma}
\newtheorem{e-proposition}[theorem]{Proposition}
\newtheorem{corollary}[theorem]{Corollary}
\newtheorem{e-definition}[theorem]{Definition\rm}

\newcommand{\BPL}{\smallskip \noindent \textbf{Proof of Lemma} }

\newcommand{\BPCOR}{\smallskip \noindent \textbf{Proof of Corollary} }

\newcommand{\BPT}{\smallskip \noindent \textbf{Proof of Theorem} }

\newtheorem{theoreme}{Th\'eor\`eme}[section]

\newtheorem{proposition}[theoreme]{Proposition}

\renewcommand{\theequation}{\arabic{equation}}
\def\og{\leavevmode\raise.3ex\hbox{$\scriptscriptstyle\langle\!\langle$~}}
\def\fg{\leavevmode\raise.3ex\hbox{~$\!\scriptscriptstyle\,\rangle\!\rangle$}}

\newcommand{\mC}{\mathrm{C}}

\newcommand{\mO}{\mathcal{O}}

\newcommand{\mM}{\mathcal{M}}

\newcommand{\mD}{\mathcal{D}}

\newcommand{\E}{\mathrm{E}}
\newcommand{\D}{\mathrm{D}}
\renewcommand{\H}{\mathrm{H}}

\newcommand{\R}{\mathbb{R}}

\newcommand{\N}{\mathbb{N}}

\newcommand{\noi}{\noindent}

\newcommand{\al}{\alpha}
\newcommand{\be}{\beta}

\newcommand{\de}{\delta}
\newcommand{\De}{\Delta}
\renewcommand{\e}{\varepsilon}

\newcommand{\Om}{\Omega}

\newcommand{\larrow}{\longrightarrow}
\newcommand{\ot}{\otimes}

\newcommand{\p}{\partial}
\newcommand{\sub}{\subseteq}
\newcommand{\set}{\setminus}
\newcommand{\by}{\times}

\newcommand{\dd}{\mathrm{d}}

\newcommand{\sgn}{\mathrm{sgn}}
\newcommand{\ess}{\mathrm{ess}}

\renewcommand{\div}{\mathrm{div}}

\def\Xint#1{\mathchoice
{\XXint\displaystyle\textstyle{#1}}%
{\XXint\textstyle\scriptstyle{#1}}%
{\XXint\scriptstyle\scriptscriptstyle{#1}}%
{\XXint\scriptscriptstyle\scriptscriptstyle{#1}}%
\!\int}
\def\XXint#1#2#3{{\setbox0=\hbox{$#1{#2#3}{\int}$ }
\vcenter{\hbox{$#2#3$ }}\kern-.6\wd0}}

\def\dashint{\Xint-}

\begin{document}
\begin{frontmatter}


\selectlanguage{english}
\title{Counterexamples in Calculus of Variations in $L^\infty$ through the vectorial Eikonal equation}


\selectlanguage{english}
\author[authorlabel1,authorlabel2]{Nikos Katzourakis},
\ead{n.katzourakis@reading.ac.uk}
\author[authorlabel3]{Giles Shaw}
\ead{giles.shaw@gmail.com}

\address[authorlabel1]{Department of Mathematics and Statistics, University of Reading, Whiteknights, PO Box 220, Reading RG6 6AX, Berkshire, England, UNITED KINGDOM}

 \thanks[authorlabel2]{N.K. has been partially financially supported through the EPSRC grant EP/N017412/1.}
\thanks[authorlabel3]{G.S. acknowledges the full financial support of the EPSRC grant EP/N017412/1.}



\begin{abstract} We show that for any regular bounded domain $\Omega\sub  \R^n$, $n=2,3$, there exist infinitely many global diffeomorphisms equal to the identity on $\p\Omega$ which solve the Eikonal equation. We also provide explicit examples of such maps on annular domains. This implies that the $\infty$-Laplace system arising in vectorial Calculus of Variations in $L^\infty$ does not suffice to characterise either limits of $p$-Harmonic maps as $p\to \infty$, or absolute minimisers in the sense of Aronsson.


\vskip 0.5\baselineskip

\selectlanguage{francais}
\noindent{\bf R\'esum\'e} \vskip 0.5\baselineskip  Nous montrons que pour tout domaine born\'e r\'egulier $ \Omega \sub \R^n $, $n=2,3$, il existe une infinit\'e de diff\'eomorphismes globaux solutions de l'\'equation iconale, \'egaux \`{a} l'identit\'e sur $ \p\Omega $. Nous donnons \'egalement des exemples explicites de telles cartes dans des domaines annulaires. Ceci implique que le syst\'eme du type $ \infty $-Laplacien apparaissant dans le Calcul des Variations vectoriel dans $ L^\infty $ ne suffit pas \`{a} caract\'eriser les limites pour $p\to \infty$ des cartes $ p $-harmoniques, ni les minimiseurs absolus au sens d'Aronsson.

\noi {\bf Contre-exemples dans le Calcul des Variations dans $ L^\infty $ par l'\'equation iconale vectorielle}

\end{abstract}
\end{frontmatter}


\selectlanguage{english}

\section{Introduction} \label{section1}

Calculus of Variations in $L^\infty$ is concerned with the variational study of supremal functionals, as well as with the necessary conditions governing their extrema. The archetypal model of interest is the functional
\begin{equation}
\label{eq:minproblem}
\E_\infty(u,\mO)\,:=\, {\ess\sup}_\mO |\D u|,\ \ \text{ for } u\in\mathrm{W}^{1,\infty}(\Omega;\R^N), \ \mO \sub \Om \text{ measurable},
\end{equation}
where $n,N\in\N$, $\Om \sub \R^n$ is a fixed open set and $\D u(x)=(\D_i u_\al(x))_{i=1...n}^{\al=1...N} \in \R^{N\times n}$ is the gradient matrix. We note that our general notation is either self-explanatory or standard. In \eqref{eq:minproblem} and throughout the remainder of this note, all norms appearing will be the Euclidean ones. In particular, on $\R^{N\times n}$ we use the norm induced by the inner product $A:B:=\smash{\sum_{\al, i}A_{\al i}B_{\al i}}$. Aronsson was the first to consider such problems in the 1960s \cite{Aronss67EFL,A4}, in the scalar-valued case $N=1$. In the general case of \eqref{eq:minproblem}, the PDE system which arises from \eqref{eq:minproblem} as the analogue of the Euler-Lagrange equations is the $\infty$-Laplace system
\begin{equation}
\label{eq:infinitylaplace}
\Delta_\infty u \, :=\,  \Big(\D u \ot \D u \,+\, |\D u|^2[\![\D u ]\!]^\bot \! \ot \mathrm I\Big) :\D^2u \, =\, 0,
\end{equation}
and has its origins in the paper \cite{Katzou12LVP}. Here, for any linear map $A :\R^n \larrow \R^N$, $[\![A]\!]^\bot$ denotes the orthogonal projection on the orthogonal complement of the range $(\mathrm R(A))^\bot \sub \R^N$. In index form, \eqref{eq:infinitylaplace} reads
\[
\sum_{1\leq \be \leq N}\sum_{1\leq i,j \leq n} \Big(\D_i u_\al \, \D_j u_\be \,+\, |\D u|^2[\![\D u ]\!]_{\al \be}^\bot \de_{ij}\Big) \D_{ij}^2u_\be \, =\, 0,\ \ \ \ 1\leq \al \leq N.
\]
A fundamental difficulty in the variational study of \eqref{eq:minproblem} is that the usual global minimisers in the space $g+\mathrm{W}_0^{1,\infty}(\Omega;\R^N)$ are not truly optimal and may not solve any PDE. To this end, the notion of absolute minimisers has been introduced. Indeed, in the scalar case it is known that absolute minimisers of \eqref{eq:minproblem} correspond uniquely to (viscosity) solutions of the \emph{scalar} version of \eqref{eq:infinitylaplace}, which reduces to $\sum_{i,j}\D_i u \D_j u \D^2_{ij}u=0$ (see e.g.\ \cite{BhDiMa89LPRE,Jensen93ULEM,ACJ,C,Kbook1}). The ``localised" concept of absolute minimisers is what forces to define \eqref{eq:minproblem} on subsets of $\Om$. In the vectorial case, the situation is more delicate and not fully understood yet, particularly when $\D u$ has rank greater than two \cite{BJW1,AyK,AK,K9}. A by now standard mechanism to study \eqref{eq:minproblem}-\eqref{eq:infinitylaplace} is through approximation by the respective $L^p$ variational notions as $p\to \infty$, namely by using
\begin{equation}
\label{LpObjects}
\E_p(u):= \|\D u\|_{L^p(\Om)},\ \text{ for } u \in\mathrm{W}^{1,p}(\Omega;\R^N) \ \ \text{ and } \ \ \ \De_pu:= \div\big(|\D u|^{p-2}\D u \big)=0,
\end{equation}
which are known as the $p$-Dirichlet functional and the $p$-Laplacian. Hence, the identification of necessary and sufficient conditions for a mapping $u\in\mathrm{W}^{1,\infty}(\Omega;\R^N)$ to occur as a (weak) limit $u$ of $p$-harmonic maps $u_p$ is of interest (see \cite{BJW1,BhDiMa89LPRE,Jensen93ULEM,K9}). Intuitively, we expect such limits to  be ``optimal'' solutions, possibly absolute minimisers of \eqref{eq:minproblem}. In the case $N=1$, a complete picture is known: the family $(u_p)_{p\geq 1}$ converges to a unique limit which is an absolute minimiser. Additionally, it follows from the form of the PDE that differentiable Eikonal functions solving $|\D u|=\mathrm{const}$, also satisfy \eqref{eq:infinitylaplace} and therefore is a $p$-harmonic limit.

On the other hand, in the case $N\geq2$, one can show the existence of infinitely-many (appropriately defined) {\it generalised $\mathrm{W}^{1,\infty}$ solutions} to \eqref{eq:infinitylaplace} which are {\it not minimising for \eqref{eq:minproblem}}, let alone absolutely minimising, see \cite{K7,CKP}. A natural question is whether this phenomenon is a defect of the notion of solution used. The principal results of this note are Theorem~\ref{thm:mainthm} and Corollary \ref{cor:mainthm}, which answer this to the negative. Accordingly, we show  for $n=N\in\{2,3\}$ the existence of infinitely many {\it arbitrarily regular} orientation preserving Eikonal diffeomorphisms $u : \Omega\larrow \R^n$ with given affine boundary conditions. These maps are {\it a fortiori $\infty$-Harmonic}, since $[\![\D u]\!]^\perp\! = 0$ when $\det(\D u)\neq 0$ and \eqref{eq:infinitylaplace} can be recast as the two independent systems
\[
\D u \,\D\big(|\D u|^2\big) = 0\ \ \text{ and } \ \  |\D u|^2[\![\D u]\!]^\perp\Delta u  = 0.
\] 
\begin{theorem}\label{thm:mainthm}
Let $\Omega\sub \R^n$ be a bounded connected domain such that  $n\in\{2,3\}$ and $\p\Omega$ is $\mC^{m+4 }$ for $m \geq2$. Then, there exist infinitely many maps $u\in \mC^m\big(\overline{\Omega};\R^n\big)$ satisfying 
\[
\text{$|\D u |\equiv \mathrm{const}$ in $\overline{\Om}$, \  $\det(\D u)>0$ in $\overline{\Om}$ \ and \ $u=\mathrm{id}$ on $\p\Omega$.}
\]
Any such $u$ is an Eikonal orientation preserving diffeomorphism, equal to the identity on the boundary.
\end{theorem}

\smallskip

\begin{corollary}
\label{cor:mainthm} Let $n,m,\Om $ be as in Theorem \ref{thm:mainthm}. Then, the Dirichlet problem for the $\infty$-Laplacian 
\[
\text{$\De_\infty u =0 $ in $\Om$ \ \ and \ \ $u=\mathrm{id}$ on $\p\Omega$,}
\]
possesses infinitely-many classical solutions $u\in \mC^m\big(\overline{\Omega};\R^n\big) \set\{\mathrm{id}\}$. In addition, none of these solutions minimises $\E_\infty(\cdot,\Om)$ among all maps in $W^{1,\infty}(\Om,\R^n)$ with $u=\mathrm{id}$ on $\p\Om$.
\end{corollary}

\smallskip

We note that the results above improve and supersede one of the main results in \cite{K6} which required $\Om$ to be a punctured ball. Since the unique solution to the Dirichlet problem for $\Delta_p u = 0$ in $\Om$ with $u=\mathrm{id}$ on $\p\Om$ is $u(x)\equiv x$ when $p<\infty$, it follows that none of our diffeomorphisms is a limit of $p$-harmonic maps as $p\to \infty$. Thus, we confirm that \eqref{eq:infinitylaplace} by itself cannot suffice to identify limits of $p$-harmonic maps and that {\it additional selection criteria are needed} to have a situation analogous to the scalar case.  

The proof of Theorem~\ref{thm:mainthm} is based on the next result of independent interest.

\smallskip

\begin{proposition}\label{prop:functionalroots}
Let $n,m,\Om$ be as in Theorem \ref{thm:mainthm}. Then, the nonlinear problem
\[
|\D u|^2  +\, 2\, \div\, u \,\equiv \,C \, \text{ in }\Om\ \ \ \text{ and } \ \ \ u=0 \,  \text{ on }\p\Om,
\]
has infinitely many non-trivial solutions $(u,C) \in (\mC^m\cap\mC^0_0)\big(\overline{\Omega};\R^n\big) \by (0,\infty)$. Additionally, the set of all solutions has the trivial solution $(0,0)$ as an accumulation point with respect to the topology of $\smash{\mC^m\big(\overline{\Omega};\R^n\big)}$.
\end{proposition}

\smallskip

Since the proofs of the above results are non constructive, we include in Section \ref{ex:counterexample} explicit examples of smooth $\infty$-Harmonic maps defined on annular domains which coincide with affine maps on the boundary.

\!\!\!\!\!\!\!\!\!\!

\section{Proofs}
\!\!\!\!

We begin with the proof of Proposition \ref{prop:functionalroots}, which is an immediate consequence of the next lemma and of the Morrey estimate, in the form of inclusion of spaces $\H^{m+2}(\Omega;\R^n)\sub \mC^m\big(\overline{\Omega};\R^n\big)$ (since $n\in\{2,3\}$).

\smallskip

\begin{lemma}\label{lem:functionalroots}
Let $n,m,\Om,$ be as in Theorem \ref{thm:mainthm} and let us define the nonlinear mapping
\[
\mM \ : \ \ (\H^{m+2}\cap\H^1_0)(\Omega;\R^n)\larrow \H^{m+1}_\sharp(\Omega)\,:= \,\left\{w\in\H^{m+1}(\Omega) :\ \int_\Omega w(x)\;\dd x = 0\right\}
\]
by setting (here the slashed integral denotes the average)
\[
\mM[u] \, := \, \frac{1}{2}|\D u|^2 \, +\, \div\, u \, -\, \frac{1}{2} \, \dashint_\Omega  |\D u(x) |^2\;\dd x.
\]
Then, the inverse image $\mM^{-1}[\{0\}]$ contains infinitely-many elements accumulating at zero. In addition, for any $\varepsilon>0$, there exists $\varphi_\e\in(\H^{m+2}\cap\H^1_0)(\Omega;\R^n)\setminus\{0\}$ such that $\mM[\varphi_\e]=0$ and $\|\varphi_\e \|_{\H^{m+2}(\Om)}<\varepsilon$.
\end{lemma}

\BPL \ref{lem:functionalroots}. First note that $\mM$ is well defined, namely its image lies in the subspace $\smash{\H^{m+1}_\sharp(\Omega)}$ of zero average. Indeed, for any $u\in (\H^{m+2}\cap\H^1_0)(\Omega;\R^n)$, the divergence theorem gives
\[
\int_\Om \mM[u](x)\,\dd x\, =\, \int_\Omega \div \, u (x)\,\dd x \, =\, \int_{\p\Omega} u(x)\cdot n(x)\,\dd\mathcal{H}^{n-1}(x) \,=\, 0,
\]
where $n : \p\Omega\larrow \R^n$ denotes the outward pointing normal vector to $\p\Omega$ and $\mathcal{H}^{n-1}$ is the $(n-1)$-Hausdorff measure. Additionally, we need to confirm that $|\D^p(|\D u|^2)| \in \mathrm L^2(\Om)$ for all $p\in\{0,...,m+1\}$. Indeed, by the Leibniz formula we have $|\D^p(|\D u|^2)|\leq \sum_{i=0}^{p}C_{i,p}|\D^{p+1-i}u| | \D^{1+i}u|$, where $C_{i,p}$ is the binomial coefficient. Since $\min\{p+1-i,1+i\}\leq m$ for all $p$, H\"older's inequality gives $|\D^{p+1-i}u| | \D^{1+i}u| \in \mathrm L^2(\Om)$ for any $i$ and $p$, because by the Sobolev inequality we have $u,|\D u|,...,|\D^m u| \in \mathrm L^\infty(\Om)$. Next, note that $\mM$ is Fr\'echet differentiable at each $u\in\smash{(\H^{m+2}\cap\H^1_0)(\Omega;\R^n)}$ with
\[
\ \ \ \ \ \mM'[u]\varphi \, =\, \D u :\D\varphi \, +\, \div \, \varphi \, -\, \dashint_\Omega (\D u :\D\varphi )(x)\;\dd x,\qquad\text{ for all }\varphi\in (\H^{m+2}\cap\H^1_0)(\Omega;\R^n).
\]
In particular, $\mM'[0]=\div$ and also $\mM'[0]$ is a bounded linear surjection from $(\H^{m+2}\cap\H^1_0)(\Omega;\R^n)$ into $\smash{\H^{m+1}_\sharp(\Omega)}$; the surjectivity of $\mM'[0]$ is a consequence of Lemma \ref{Lemma2.2} that follows. Next, since $\ker(\mM'[0])=\big\{v\in(\H^{m+2}\cap\H^1_0)(\Omega;\R^n) : \div \, v\equiv 0 \big\}$ is a closed subspace of the Hilbert space $(\H^{m+2}\cap\H^1_0)(\Omega;\R^n)$, it possesses an orthogonal complement $V\sub (\H^{m+2}\cap\H^1_0)(\Omega;\R^n)$:
\[
 \ker(\mM'[0]) \oplus V \, =\, (\H^{m+2}\cap\H^1_0)(\Omega;\R^n).
\]
By noting that $\mM'[0]|_{ V} : V\larrow  \smash{\H^{m+1}_\sharp(\Omega)}$ is a linear isomorphism, the canonical isomorphism between $\ker(\mM'[0])\oplus V$ and $\ker(\mM'[0])\times V$ allows us to view $\mM$ as a map on $\ker(\mM'[0])\times V$ by setting $\mM[(u,v)]:=\mM[u+v]$. Then, the implicit function theorem (see e.g.\ \cite[Th.\ 4.E]{Zeidle95AFA}) implies that, for $\e>0$ small enough, there exists a continuous map $\gamma :  \ker(\mM'[0]) \cap \{v: \|v \|_{\H^{m+2}(\Om)}<\e\} \larrow V$ with $\gamma(0)=0$ and
\[
\mM[\varphi+\gamma(\varphi)] = 0 \quad\text{ for all }\varphi\in\ker(\mM'[0])\text{ with }\|\varphi\|_{\H^{m+2}(\Om)}<\e.
\]
Consequently, since $\ker(\mM'[0]) \neq \{0\}$ (as for instance $\mathrm{curl}^*\psi \in \ker(\mM'[0])$ for any $\psi\in\mC^\infty_c(\Omega;\R^{n\times n}_{\mathrm{skew}})$) and $\gamma$ is continuous with $\gamma(0)=0$ we deduce that, for every $\e>0$, there exists $\varphi_\e\in(\H^{m+2}\cap\H^1_0)(\Omega;\R^n)$ such that $\mM[\varphi_\e]=0$ and $\|\varphi_\e\|_{\H^{m+2}(\Om)}<\e$.
\qed 
\smallskip

The next result completes the proof of Lemma \ref{lem:functionalroots}.
\smallskip

\begin{lemma} \label{Lemma2.2} For any $f\in \H^{m+1}_\sharp(\Omega)$, the next Dirichlet problem admits a solution in $(\H^{m+2}\cap\H^1_0)(\Omega;\R^n)$:
\[
\div\, u = f \ \text{ in }\Omega \ \ \text{ and } \ \ u =0 \ \text{ on }\p\Omega.
\] 
\end{lemma}

\BPL \ref{Lemma2.2}. The claim follows from the Sobolev version of arguments presented in~\cite[Ch.\ 9]{CsDaKn12TPB} and standard regularity results for the Neumann problem, which we sketch briefly for completeness. Since $\p\Omega$  is assumed to be $\mC^{m+4 }$, regularity theory for Poisson's equation implies that we can always find $w\in \H^{m+3}(\Omega)$ such that $\Delta w=f$ in $\Omega$ and $\D_n w = 0$ on $\p\Omega$. It remains to show that we can find $b\in\H^{m+2}(\Omega;\R^{n})$ such that $\div \, b= 0$ in $\Omega$ and $b=-\D  w$ on $\p\Omega$ (since we can then take $u := \D  w -b$ as our desired solution). To this end, let us fix $(i,j)\in\{1,\ldots,n\}^2$ and set $c_{ij}:=(n_j\, \D_i w -n_i\, \D_j w)n$. Consider then the Biharmonic function $d_{ij}\in \H^{m+3}(\Omega)$ solving the Dirichlet problem (see \cite[Th.\ 2.2]{GGS})
\[
\Delta^2 d_{ij} =0 \ \text{in }\Omega \ \ \text{ and }\ \ d_{ij}=0\ \text{on }\p\Omega \ \ \text{ and }\ \ \D  d_{ij} = c_{ij}\ \text{on }\p\Omega.
\]
Defining $b:=\mathrm{curl}^*d \in\H^{m+2}(\Omega;\R^n)$, where $\big(\mathrm{curl}^*d\big)_i := \sum_{j<i}\D_j  d_{ji}-\sum_{j>i}\D_j d_{ij}$, we see that $\div\, b =0$ in $\Omega$ and by using that $\D_n w=0$ on $\p\Omega$, we can easily confirm that $b=-\D  w$ on $\p\Omega$. 
\qed

\smallskip

Now we may establish our main result.

\BPT \ref{thm:mainthm}. By continuity of the determinant, there exists $\varepsilon>0$ such that $|A|<\varepsilon$ implies $\det(\mathrm I +A)>\frac{1}{2}$. Using Proposition~\ref{prop:functionalroots}, we can find $\varphi\in\mC^m\big(\overline{\Omega};\R^n\big)\set\{0\}$ with $\phi|_{ \p\Omega}=0$ and $\|\D\varphi\|_{C^0(\Om)}<\e$, satisfying
\[
|\D\varphi |^2 +2\div \varphi \, \equiv \, C\ \ \text{ and }  \ \
\det(\mathrm{I} +\D\varphi )  > {1}/{2}\ \text{ in }\Omega,
\]
for some $C>0$. Defining $u\in\mC^m\big(\overline{\Omega};\R^n\big)$ by $u:=\mathrm{id} + \varphi$, we have  $u=\mathrm{id}$ on $\p\Omega$. Additionally,
\begin{align*}
|\D u |^2=|\mathrm I  + \D\varphi |^2 \, &=\, |\mathrm I|^2\,+\,2\,\mathrm I\!:\!\D\varphi \, +\, |\D\varphi |^2\, =\, n^2 \, +\, \big(2\,\div \,\varphi \, +\, |\D\varphi |^2\big)\, =\, n^2 +C
\end{align*}
and also $\det(\D u)\geq 1/2$ on $\Om$, as required. Evidently, $u$ is a local diffeomorphism from $\Om$ into $\R^n$. The fact that $u$ is a global diffeomorphism follows from standard degree theory results (see e.g.\ \cite[Th.\ 19.12]{CsDaKn12TPB}).
\qed

\BPCOR \ref{cor:mainthm}. Evidently, for any $u$ as above we have $\De_\infty u=0$ in $\Om$ and $u=\mathrm{id}$ on $\p\Om$. Further, since $|\D u |^2\equiv n^2 +C= |\mathrm I|^2+C>|\mathrm I|^2$ and also $\D (\mathrm{id})=\mathrm I$, we obtain $\E_\infty(u,\Om)>\E_\infty(\mathrm{id},\Om)$. \qed

\!\!\!\!\!\!\!\!

\section{Explicit constructions}
\label{ex:counterexample}

\!\!\!\!\!

\begin{lemma} \label{explicit}

Let $n\in2\N$, $\Om :=\{x\in\R^{n}: 1<|x|< e^{2\pi}\}$ and $S\in\R^{n\times n}$ an orthogonal, skew-symmetric matrix whose spectrum satisfies $\sigma(S)\sub \{\pm i, 0\}$ so that $e^{2\pi S} =\mathrm I$. Let $u : \overline{\Om}\larrow \R^{n}$ be given by $u(x):=\mathrm{e}^{\log(|x|)S}x$. Then, $u \in \mC^\infty\big(\overline{\Om};\R^{n}\big)\set\{\mathrm{id}\}$, $u=\mathrm{id}$ on $\p\Om$, $|\D u|^2\equiv n^2+1$ in $\overline{\Om}$ and $\det(\D u)\equiv 1$ in $\overline{\Om}$. In particular, $u$ is a global $\infty$-Harmonic orientation preserving diffeomorphism.
\end{lemma}

\BPL \ref{explicit}. It is clear that $u \in \mC^\infty\big(\overline{\Om};\R^{n}\big)\set\{\mathrm{id}\}$, $u=\mathrm{id}$ on $\p\Om$. By using standard properties of the matrix exponential (in particular that $e^{f(t)S}S=Se^{f(t)S}$ and $\D_t(e^{f(t)X})=f'(t)Se^{f(t)S}$ for any $f\in \mC^1(\R;\R)$) and setting for convenience $\sgn(x):=x/|x|$, when $x\in\R^n \set\{0\}$, we easily compute that
\[
\D u(x)  = \mathrm{e}^{\log(|x|)S}\big(\mathrm I \, +\, (S \, \sgn(x)) \otimes \sgn(x)\big), \ \ \text{ for }x\in\Om.
\]
Since $\mathrm{e}^{\log(|x|)S}$ is orthogonal and $|OA|=|A|$ for any $A,O\in\R^{n\times n}$ with $O$ being orthogonal, we have
\begin{align*}
|\D u(x) |^2 \, = \, \big|\mathrm I + (S \, \sgn(x))\otimes \sgn(x) \big|^2\, =\, n^2 \,+ \, 2\left(S\, \sgn(x)\right) \cdot \sgn(x)  \,+\,\left|S \,\sgn(x)\right|^2\left|\sgn(x)\right|^2.
\end{align*}
Because $S$ is both skew-symmetric and orthogonal, we have $(Se)\cdot e=0$ and also $|Se|=1$ when $|e|=1$. We therefore have $|\D u|^2 \equiv n^2 + 1$ on $\Om $. By using once again that $\mathrm{e}^{\log(|x|)S}$ is orthogonal, we have
\[
\det(\D u(x)) \, =\, \det\big(\mathrm{e}^{\log(|x|)S}\big)\det \!\big(\mathrm{I}  + (S\, \sgn(x))\otimes \sgn(x)\big) \, =\, \det\!\big(\mathrm{I}  + (S\, \sgn(x))\otimes \sgn(x)\big).
\]
By the Matrix Determinant Lemma, $\det(\mathrm{I}+a\otimes b)=1+a\cdot b$ for any $a,b\in\R^n$ and so we can use again the skew-symmetry of $S$ to deduce $\det(\D u(x)) =1+ (S\sgn(x))\cdot \sgn(x) = 1$ for any $x\in\Om$.  
\qed

\smallskip

\begin{lemma} \label{lemma} For any $n,N \geq 2$, there exists an explicit smooth $\infty$-Harmonic map defined on a cylindrical subdomain of $\R^n$ with values in $\R^N$ which coincides with an affine map on the boundary of the domain.
\end{lemma}

\BPL \ref{lemma}. Let $u : \R^2 \supseteq \overline{\Om} \larrow \R^2$ be the mapping constructed in Lemma \ref{explicit}. Then, by setting $v(x):=(u(x),0)^\top$, we obtain a  map $v: \Om \larrow \R^{2+k}$ for any $k\in\N$ with the desired properties. Indeed, we have $|\D v|^2=|\D u|^2$, which gives $\D v \ot \D v :\D^2 v \equiv 0$. Further, for any $x\in\Om$ we have $\mathrm R(\D v(x))=\R^2 \by \{0\}$ and $\De v(x) \in \R^2 \by \{0\}$, which gives $[\![\D v]\!]^\bot\De v\equiv 0$. Hence, $\De_\infty v \equiv 0$ in $\Om$, whilst $v=(\mathrm{id},0)$ on $\p\Om$. 

Further, by setting $w(x,y):=v(x)$, we obtain a  map $w: \R^{2+l} \supseteq \Om\by\R^l \larrow \R^{2+k}$ for any $k,l\in\N$ defined on a cylindrical annulus with $|\D w|^2=|\D v|^2$, which gives $\D w \ot \D w :\D^2 w \equiv 0$. Also, for any $(x,y)\in\Om\by \R^l$ we have $\mathrm R(\D w(x,y))=\mathrm R(\D v(x))$ and $\De w(x,y)=\De v(x)$, giving $[\![\D w]\!]^\bot\De w\equiv 0$ and thus $\De_\infty w \equiv 0$ in $\Om \by \R^l$. Finally, note that $w=(\mathrm{Proj}_{\R^2},0)$ on $\p(\Om \by \R^l)$. 
\qed

\smallskip

\noi \textbf{Acknowledgements.} N.K. would like to thank Roger Moser for numerous inspiring scientific discussion on Calculus of Variations in $L^\infty$. Additionally, both authors have benefited whilst preparing this note by the ideas emerging in an unpublished existential counterexample relevant to our results herein, which was kindly and selflessly shared with the authors. We are grateful to Roger Moser for this deep insight. We would also like to thank the referee for their constructive comments and suggestion which improved the content and the presentation. In particular, we are grateful for a correction in the proof of Lemma 2.1.

\!\!\!\!\!

\end{document}